\documentclass[12pt]{article}
\usepackage{amsmath}
\usepackage{amssymb}
\usepackage{amsthm}
\usepackage{epsfig}

\begin{document}
\theoremstyle{plain}
\newtheorem{theorem}{Theorem}
\newtheorem{prop}[theorem]{Proposition}
\newtheorem{corollary}[theorem]{Corollary}
\newtheorem{lemma}[theorem]{Lemma}
\newtheorem{question}[theorem]{Question}
\newtheorem{conjecture}[theorem]{Conjecture}
\newtheorem{assumption}[theorem]{Assumption}

\theoremstyle{definition}
\newtheorem{definition}[theorem]{Definition}
\newtheorem{notation}[theorem]{Notation}
\newtheorem{condition}[theorem]{Condition}
\newtheorem{example}[theorem]{Example}
\newtheorem{introduction}[theorem]{Introduction}

\numberwithin{theorem}{section}

\makeatletter                          
\let\c@equation\c@theorem              
\makeatother                           
\renewcommand{\theequation}{\arabic{section}.\arabic{equation}}

\makeatletter                          
\let\c@figure\c@theorem              
\makeatother                           
\renewcommand{\thefigure}{\arabic{section}.\arabic{figure}}

\newcommand{\todo}[1]{\vspace{5 mm}\par \noindent
   \marginpar{\textsc{ToDo}}\framebox{\begin{minipage}[c]{0.95 \textwidth}
   \tt #1 \end{minipage}}\vspace{5 mm}\par}

\providecommand{\abs}[1]{\lvert#1\rvert} \providecommand{\norm}[1]{\lVert#1\rVert}
\providecommand{\Z}{\mathbb{Z}} \providecommand{\R}{\mathbb{R}}
\providecommand{\N}{\mathbb{N}} \providecommand{\C}{\mathbb{C}}
\providecommand{\Q}{{\mathbb{Q}}} \providecommand{\x}{\mathbf{x}}
\providecommand{\y}{\mathbf{y}} \providecommand{\z}{\mathbf{z}}
\long\def\symbolfootnote[#1]#2{\begingroup%
\def\thefootnote{\fnsymbol{footnote}}\footnote[#1]{#2}\endgroup}

\title{Neighborhood complexes and generating functions for affine semigroups}

\author{Herbert E. Scarf and Kevin M. Woods\footnote{This author was
partially supported by an NSF Graduate Research Fellowship.}}

\maketitle

\begin{abstract}
Given $a_1,a_2,\ldots,a_n\in\Z^d$, we examine the set, G, of all non-negative
integer combinations of these $a_i$.  In particular, we examine the generating
function $f(\z)=\sum_{b\in G}\z^b$.  We prove that one can write this generating
function as a rational function using the neighborhood complex (sometimes called the
complex of maximal lattice-free bodies or the Scarf complex) on a particular lattice
in $\Z^n$. In the generic case, this follows from algebraic results of D. Bayer and
B. Sturmfels. Here we prove it geometrically in all cases, and we examine a
generalization involving the neighborhood complex on an arbitrary lattice.
\end{abstract}

\section{Introduction}
Given positive integers $a_1,a_2,\ldots,a_n$, let
\[G=\{\mu_1a_1+\mu_2a_2+\cdots+\mu_na_n:\  \mu_i\in\Z_{\ge 0}\}.\]
In other words, $G$ is the additive semigroup (with zero) generated by
$a_1,a_2,\ldots,a_n$.  If the greatest common divisor of $a_1, a_2,\ldots,a_n$ is
one, then all sufficiently large integers are in $G$, and the Frobenius problem is
to find the largest integer not in $G$.  We would like to say something about the
structure of the set $G$.  In particular, define the generating function
\[f(z)=\sum_{b\in G}z^b.\]
This generating function converges for $\abs{z}<1$.  We would like to calculate
$f(z)$ in a nice form.  It will turn out that we can obtain it from the neighborhood
complex (sometimes called the Scarf complex or the complex of maximal lattice-free
bodies; we will define it shortly) of an associated lattice. This was proved by D.
Bayer and B. Sturmfels using algebraic methods \cite{BS98}. Here we prove it
geometrically.

We do not need to restrict ourselves to the case where $G$ is one dimensional. In
general, Let $A$ be a $d \times n$ matrix of integers with columns
$a_1,a_2,\ldots,a_n \in \mathbb{Z}^d$, and define
\[G=\{\mu_1a_1+\mu_2a_2+\cdots+\mu_na_n:\  \mu_i\in\Z_{\ge 0}\}=\{A\xi:\ \xi\in\Z^n_{\ge 0}\}.\]
Then the $d=1$ case corresponds to the Frobenius problem.  Define the generating
function
\[f(\mathbf{z})=\sum_{b\in G} \mathbf{z}^b=\sum_{b=(b_1,\ldots,b_d)\in G}z_1^{b_1}z_2^{b_2}\cdots z_d^{b_d},\]
where $\mathbf{z}=(z_1,z_2,\ldots,z_d)$.  We assume that there exists an
$l=(l_1,l_2,\ldots,l_d)\in \mathbb{R}^d$ such that $\langle l,a_i\rangle<0$ for all
$i$.  Then for all $\mathbf{z}$ in a neighborhood of
$(e^{l_1},e^{l_2},\ldots,e^{l_d})$ we have $\norm{\mathbf{z}^{a_i}}<1$, and so
$f(\mathbf{z})$ will converge in this neighborhood.  Note that if there were no such
$l$, then $G$ would contain a linear subgroup, and $f(\z)$ would not converge on any
open subset of $\C^d$.  Since the structure of a linear group is simple, however, we
are not concerned with such $G$.

We would like to calculate this generating function, $f(\z)$. Theorem \ref{GF} gives
the answer.

Let $\Lambda \subset \mathbb{Z}^n$ be the lattice
\[\{\lambda\in \mathbb{Z}^n :\  A \lambda=0\}.\]
We will shortly define the \emph{neighborhood complex}, $S$, a simplicial complex
whose vertices are $\Lambda$. By a simplicial complex, we mean that $S$ is a
collection of finite subsets of $\Lambda$, and that if $s\in S$, then all subsets of
$s$ are also in $S$. The vertices of the complex are the $\{s\}\in S$, the edges are
the $\{s,s'\}\in S$, and so forth.  In this paper, we will not count the empty set
as a simplex of $S$.  This complex will not, in general, be geometrically realizable
in the linear span of $\Lambda$.

For $s=\{\lambda^1,\lambda^2,\ldots,\lambda^k\}$ with $\lambda^i \in \Lambda$,
define
\[\max(s)=\max\{\lambda^1,\lambda^2,\ldots,\lambda^k\},\]
where the maximum is taken coordinate-wise (for example,
$\max\left\{(1,-1),(0,0)\right\}=(1,0)$). We say that $\Lambda$ is \emph{generic}
if, whenever some nonzero $\lambda=(\lambda_1,\lambda_2,\ldots,\lambda_n)\in\Lambda$
has $\lambda_i=0$, for some $i$, then there is a $\lambda'\in\Lambda$ with
$\lambda'<\max(\lambda,0)$.

When $\Lambda$ is generic, define $S$, as follows.  We have $s$ is in $S$ if and
only if for no $\lambda \in \Lambda$ is $\lambda<\max(s)$. If $s\in S$, then all
subsets of $s$ are in $S$ as well, so $S$ is a simplicial complex.  In Section 5, we
will define $S$ in the non-generic case. Basically, we must \emph{perturb} the
vertices slightly so we are in the generic case.

\begin{figure}
        \begin{center}
        \includegraphics[height= 3in]{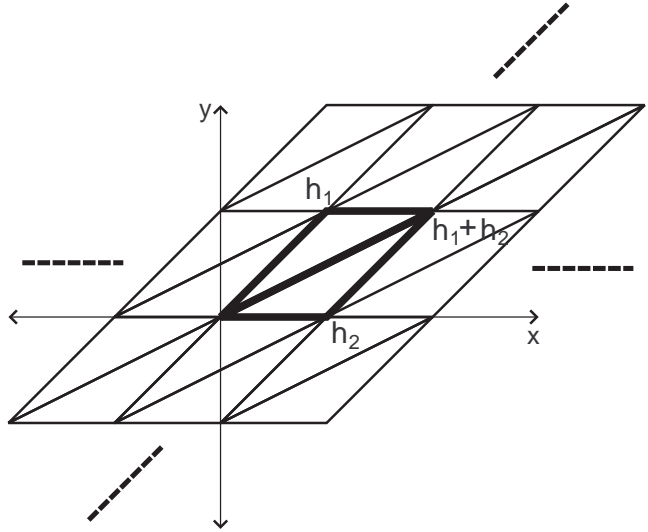}
        \caption{Example \ref{Ex1}, the neighborhood complex (transformed to lie on $\Z^2$)
        when $A$ is a $1\times 3$ matrix.}
        \label{Fig1}
        \end{center}
\end{figure}

\begin{example}\label{Ex1} If $A$ is a $1\times 3$ matrix (that is, $G$ is a one-dimensional additive semigroup
with three generators), then $\Lambda$ is a two-dimensional sublattice of $\Z^3$.
There exists a basis $\{h^1,h^2\}$ of $\Lambda$ such that the neighborhood complex
consists of vertices $\{x\}$, for $x\in\Lambda$; edges $\{x,x+h^1\},$ $\{x,x+h^2\}$,
and $\{x,x+h^1+h^2\}$; and triangles $\{x,x+h^1,x+h^1+h^2\}$ and
$\{x,x+h^2,x+h^1+h^2\}$ (see Figure \ref{Fig1}, where $\Lambda$ has been transformed
to be $\Z^2$, and see, for example, \cite{S97}). Notice that these triangles exactly
tile the linear span of $\Lambda$. This is not true in higher dimensions.
\end{example}

Neighborhood complexes often appear in integer programming in a slightly different,
but equivalent, form.  Let $r$ be the dimension of the lattice $\Lambda$, and let
$B$ be an $n\times r$ integer matrix whose columns form a basis for $\Lambda$, so
that $\Lambda=\{Bx:\ x\in\Z^r\}$.  Then we may form a simplicial complex, $S'$, on
$\Z^r$, as follows.  Given $s'=\{x^1,x^2,\ldots,x^k\}\subset \Z^r$, let $P_{s'}$ be
the polytope defined by
\[P_{s'}=\big\{x\in\R^{r}:\ Bx\le\max\{Bx^1,Bx^2,\ldots,Bx^k\}\big\}.\]
$P_{s'}$ is the smallest polytope of any $\{x:\ Bx\le b\}$, for $b\in\Z^n$, which
contains $s'$.  In the generic case, we say that $s'\in S'$ if and only if $P_{s'}$
contains no integer points in its interior.  It is easily seen that $S'$ and $S$ are
isomorphic under the map $x\mapsto Bx$.

If $\{0,x'\}\in S'$ is an edge of the complex, then $x'$ is called a \emph{neighbor}
of the origin.  The set of neighbors of the origin form a test set for the family of
integer programs
\[\text{minimize }\langle\beta_n,x\rangle \text{ such that }\langle\beta_i,x\rangle\le b_i, \text{ for
} 1\le i\le n-1,\] where $\beta_i$ is the $i$th row of $B$ and
$b=(b_1,b_2,\ldots,b_{n-1})$ is allowed to vary in $\Z^{n-1}$, and where
$\langle\cdot,\cdot\rangle$ is the standard dot product on $\R^r$. The set of
neighbors is a test set, because, for a fixed $b$, if $x$ is a feasible solution
(that is, it satisfies the linear inequalities), then $x$ minimizes
$\langle\beta_n,x\rangle$ if and only if there is no neighbor, $x'$, of the origin
such that both $x-x'$ is feasible and
$\langle\beta_n,x-x'\rangle<\langle\beta_n,x\rangle$. For an introduction to
neighbors and their applications to integer programming, see \cite{S97}.

Returning to $S$, the complex with vertices in $\Lambda$, we see that $S$ is
invariant under translation by $\Lambda$. Let $\bar{S}$ be a set of distinct
representatives of $S$ modulo $\Lambda$.  Let
\[f_{\bar{S}}(\mathbf{z})=\frac{\sum_{s\in\bar{S}}
(-1)^{\dim(s)}\mathbf{z}^{A\cdot \max(s)}}{\prod(1-\mathbf{z}^{a_i})}.\] The
following theorem states that this is the generating function that we are looking
for.

\begin{theorem}
\label{GF} Given a $d\times n$ matrix of integers $A$, let
$\Lambda=\{\lambda\in\Z^n:\ A\lambda=0\}$.  Define the neighborhood complex $S$ on
$\Lambda$ as above (we define $S$ in the non-generic case in Section 5), and let
$\bar{S}$ be a set of distinct representatives of $S$ modulo $\Lambda$.  If
$f(\mathbf{z})$ and $f_{\bar{S}}(\mathbf{z})$ are defined as above, then
\[f(\mathbf{z})=f_{\bar{S}}(\mathbf{z}).\]
\end{theorem}

In the generic case, this theorem follows from algebraic results of D. Bayer and B.
Sturmfels \cite{BS98}, but we prove it here using elementary geometric methods.
Bayer and Sturmfels construct the \emph{hull complex}, which coincides with
$\bar{S}$ when $\Lambda$ is generic, but which is larger than $\bar{S}$ in the
non-generic case.  Note that they use Hilbert series terminology, which is
equivalent, because $f(\z)$ is the Hilbert series for the monomial ring
$\C[\x^{a_1},\x^{a_2},\ldots,\x^{a_n}]$ with the standard $\Z^d$-grading.  A.
Barvinok and K. Woods show \cite{BW03} that $f(\z)$ can be written as a ``short''
rational generating function (much shorter than $f_{\bar{S}}(\z))$, but, when
written in that form, the structure of the neighborhood complex is lost.

The function $f_{\bar{S}}(\z)$ makes sense even if $\Lambda$ is a proper sublattice
(perhaps of full dimension, perhaps not) of $\{\lambda\in \mathbb{Z}^n :\  A
\lambda=0\}$. That is, we may still define the neighborhood complex, $S$, and then
take $\bar{S}$, a set of distinct representatives of $S$ modulo $\Lambda$, and
define $f_{\bar{S}}(\z)$ as above. Does $f_{\bar{S}}$ have an interpretation as a
generating function, as in Theorem \ref{GF}?

In fact, it does, as follows.  Let $\Lambda$ be any lattice in $\Z^n$ such that $A
\lambda =0$, for all $\lambda\in\Lambda$.  Given $b\in\Z^d$, define
\[T_b=\{\xi\in\Z^n:\  \xi\ge 0 \text{ and } A\xi=b\}.\]
That is, $T_b$ represents the set of ways to write $b$ as a nonnegative integer
combination of the $a_1,a_2,\ldots,a_n$ (and so $T_b$ is nonempty if and only if $b$
is in the semigroup $G$).  Define an equivalence relation on $T_b$ by
\[\xi\sim\eta \text{ iff }\xi-\eta\in\Lambda.\]
Let $c_b$ be the number of equivalence classes in $T_b$.  Then we have the following
theorem, which says that the $c_b$ are the coefficients of the Laurent power series
$f_{\bar{S}}(\z)$.

\begin{theorem}
\label{GenGF} Given a $d\times n$ matrix of integers $A$ and a lattice
$\Lambda\in\Z^n$ such that $A\lambda=0$ for all $\lambda\in\Lambda$, define the
neighborhood complex $S$ on $\Lambda$ as above, and let $\bar{S}$ be a set of
distinct representatives of $S$ modulo $\Lambda$.  If $f_{\bar{S}}(\z)$ and $c_b$
are defined as above, then
\[f_{\bar{S}}(\z)=\sum_{b\in \Z^d}c_b\z^b.\]
\end{theorem}

When $\Lambda$ is a generic lattice, this theorem can be retrieved from a result of
I. Peeva and B. Sturmfels \cite{PS98}, but they again use algebraic methods.  In the
case where $\Lambda$ is the full lattice $\{\lambda\in \mathbb{Z}^n :\  A
\lambda=0\}$, every element of $T_b$ is equivalent to every other, since if
$A\xi=A\eta =b$, then $A\cdot(\xi-\eta)=0$ and so $\xi-\eta\in \Lambda$. In this
case, if $b\in G$ then $c_b=1$ (and if $b\notin G$ then $c_b=0$), and we recover
Theorem \ref{GF}.

At the other extreme, if $\Lambda=\{0\}$, each element of $T_b$ is in its own
equivalence class.  Then, since $S$ is the complex with one vertex $0$, we have
\[f_{\bar{S}}(\z)=\frac{1}{\prod(1-\mathbf{z}^{a_i})},\]
and Theorem \ref{GenGF}, in this case, is clear.  We will present other examples of Theorems \ref{GF} and \ref{GenGF} in Section 2.

Let $ L$ be the full orthogonal lattice $\{\ell\in \mathbb{Z}^n :\  A \ell=0\}$, and
let $C$ be a lattice invariant simplicial complex on $ L$. Note that when the
lattice $\Lambda$ in Theorem \ref{GenGF} is not all of $ L$, then $S$ itself is not
$ L$-invariant. In this case, if $\bar{L}$ is a set of distinct representatives of
$L$ modulo $\Lambda$, then the complex $C$ we will examine will be the disjoint
union
\[C=\bigcup_{\ell\in \bar{L}}S+\ell,\]
where
\[S+\ell=\big\{\{\lambda^1+\ell,\lambda^2+\ell,\ldots,
\lambda^k+\ell\}:\  \{\lambda^1,\lambda^2,\ldots,\lambda^k\}\in S\big\},\] which is
$ L$-invariant.

Define $C_\xi$ to be the subcomplex of $C$ consisting of simplices $s\in C$ such
that
\[\max(s)\le\xi.\]
$C_\xi$ is a simplicial complex, though it need not be pure (that is, its maximal
simplices may not all be of the same dimension).

\begin{figure}
        \begin{center}
        \includegraphics[height= 3in]{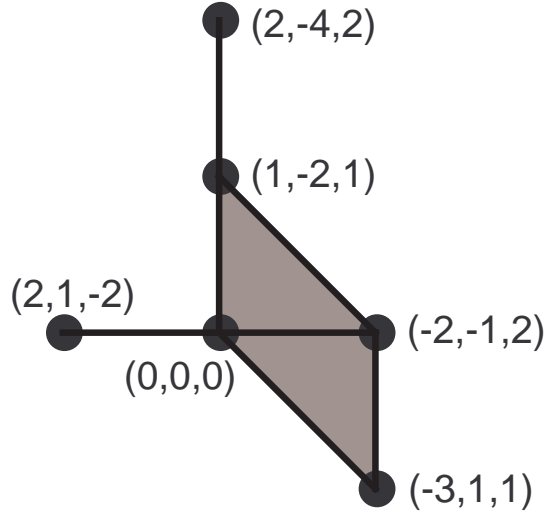}
        \caption{$C_\xi$ from Example \ref{ExCXI}}
        \label{FigCXI}
        \end{center}
\end{figure}

\begin{example}
\label{ExCXI} Let $A=[3\ 4\ 5]$ so that $G$ is the additive semigroup generated by
3, 4, and 5. Let $\Lambda=L=\{\lambda\in \mathbb{Z}^3 :\  A \lambda=0\}$, and let
$C=S$ be the neighborhood complex defined on $\Lambda$.  One can show that $C$ has
vertices $\{\lambda\}$, for $\lambda\in\Lambda$; edges
$\{\lambda,\lambda+(2,1,-2)\},$ $\{\lambda,\lambda+(1,-2,1)\}$, and
$\{\lambda,\lambda+(3,-1,-1)\}$; and triangles
$\{\lambda,\lambda+(2,1,-2),\lambda+(3,-1,-1)\}$ and
$\{\lambda,\lambda+(1,-2,1),\lambda+(3,-1,-1)\}$.

Let $\xi=(2,1,2)$.  Then Figure \ref{FigCXI} shows $C_{\xi}$.  Note that, if
$\{\lambda\}$ is a vertex of $C_{\xi}$, then $\xi-\lambda\ge 0$ by definition of
$C_\xi$, and
\[A(\xi-\lambda)=A\xi=20,\]
and so (as will be important later), each vertex of $C_{\xi}$ corresponds to a way
to write 20 as a nonnegative integer combination of 3,4, and 5.  For example,
$\{(1,-2,1)\}$ is a vertex of $C_{\xi}$, $\xi-(1,-2,1)=(1,3,1)$, and $1\cdot
3+3\cdot 4+1\cdot 5=20$.
\end{example}

Define the Euler characteristic, EC$(C_\xi)$, by
\[\text{EC}(C_\xi)=\sum_{s\in C_\xi}(-1)^{\dim(s)}.\]

Since $C$ is $ L$-invariant, $C_{\xi-\ell}=C_\xi - \ell$ for all $\ell\in
 L$. Then, given $b\in\Z^d$, all of the $C_\xi$, for $\xi\in\Z^n$ such
that $A\xi=b$, are isomorphic to each other, and we can define
\[d_b=\text{EC}(C_\xi),\text{ for some (any) }\xi\text{ such that } A\xi=b.\]

 We will prove Theorems \ref{GF}
and \ref{GenGF} using the following lemma, which says that these Euler
characteristics, $d_b$, are the coefficients of the Laurent power series
$f_{\bar{C}}(\mathbf{z})$.

\begin{lemma}
\label{GF2} If $A$ is a $d\times n$ matrix of integers and $C$ is a lattice
invariant simplicial complex on $ L=\{\ell\in \mathbb{Z}^n :\  A \ell=0\}$, let
$d_b$ be defined as above, for all $b\in\Z^d$. If $\bar{C}$ is a set of distinct
representatives of $C$ modulo $ L$, then
\[f_{\bar{C}}(\mathbf{z})=\sum_{b\in\Z^d}d_b\z^b,\]
where
\[f_{\bar{C}}(\mathbf{z})=\frac{\sum_{s\in\bar{C}}
(-1)^{\dim(s)}\mathbf{z}^{A\cdot
\max(s)}}{\prod(1-\mathbf{z}^{a_i})}.\]
\end{lemma}

 We will prove this lemma in Section 3. First, in Section 2, we will give several examples
 of Theorems \ref{GF} and \ref{GenGF}.  In Section 4, we examine neighborhood complexes and
 make the Euler characteristic calculations necessary to prove Theorems \ref{GF} and
 \ref{GenGF} from Lemma \ref{GF2}.  The key ingredient in these calculations will be the fact
 (first proved in \cite{BSS98}) that these neighborhood complexes have a very nice topological
 structure. In Section 5, we examine the non-generic case, and
 prove Theorems \ref{GF} and \ref{GenGF} for these lattices.

\section{Examples}
In this section, we look at several examples of Theorems \ref{GF} and \ref{GenGF}.
First we examine Theorem \ref{GF}, for varying $n$ and $d$.

Suppose $d=1$.  If $a_1,a_2,\ldots,a_n$ are positive integers whose greatest common
divisor is one, then the Frobenius number is the largest integer not in $G$.  The
problem of finding this number dates back to Frobenius and Sylvester.  H. Scarf and
D. Shallcross \cite{SS93} have related the Frobenius number itself to the
neighborhood complex.  They show (using slightly different terminology) that, if
\[N=\max\{A\cdot \max(s):\  s\text{ is in the neighborhood complex, } S\},\]
then the Frobenius number is
\[N-(a_1+a_2+\cdots +a_n).\]
Note that, in the terminology of this paper, $N$ is the largest exponent in the
numerator of $f_{\bar{S}}(\z)=f(\z)$.

\begin{figure}
        \begin{center}
        \includegraphics[height= 3in]{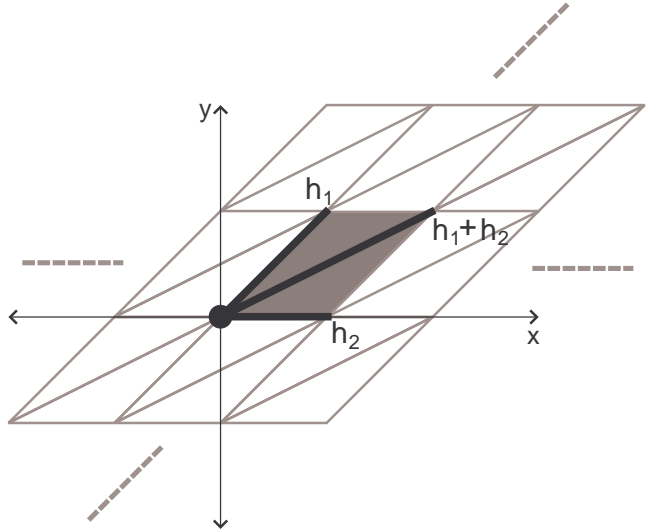}
        \caption{$\bar{C}$ (transformed to have vertices $\Z^2$), for $d=1, n=3$ in Example \ref{Ex2}}
        \label{Fig2}
        \end{center}
\end{figure}

\begin{example}
Theorem \ref{GF}, with $d=1, n=2$.  Then
\[f(z)=\frac{1-z^{\text{lcm}(a_1,a_2)}}{(1-z^{a_1})(1-z^{a_2})}.\]
\end{example}
In this case, we may choose $\bar{S}$ to consist of the vertex $\{0\}$ and the edge
$\{0,h\}$, where $h$ is a generator of the lattice $\Lambda$. This formula can
easily be verified directly.

\begin{example}
\label{Ex2} Theorem \ref{GF}, with $d=1, n=3$.  Then
\[f(z)=\frac{poly(z)}{(1-z^{a_1})(1-z^{a_2})(1-z^{a_3})},\]
where $poly(z)$ is a polynomial with at most 6 monomials.
\end{example}
In this case, $\bar{S}$ consists of one vertex, three edges, and two triangles (see
\cite{S81}). More specifically, for some $h^1,h^2\in\Lambda$, we may take $\bar{S}$
to be the set with vertex $\{0\}$; edges $\{0,h^1\}$, $\{0,h^2\}$, and
$\{0,h^1+h^2\}$; and triangles $\{0,h^1,h^1+h^2\}$ and $\{0,h^2,h^1+h^2\}$ (see
Figure \ref{Fig2}). This formula was previously shown in \cite{D96}, and also
follows from \cite{H70}, but their proofs required algebraic methods.

Here is a specific example:
\begin{example}
Theorem \ref{GF}, with $a_1=11$, $a_2=17$, and $a_3=23$.  Then
\[f(z)=\frac{1-z^{34}-z^{138}-z^{132}+z^{155}+z^{149}}{(1-z^{11})(1-z^{17})(1-z^{23})}.\]
\end{example}
In this case, we may take $h^1=(1,-2,1)$ and $h^2=(11,1,-6)$.

Unfortunately, for $d=1, n\ge 4$, the number of simplices in $\bar{S}$ may be very
large, so no formula is quite so nice. Now we examine Theorem \ref{GF} for arbitrary
$d$.

\begin{example}
Theorem \ref{GF}, with n=$d+1$.  If the $\R$-span of the $a_i$ is all of $\R^d$ (and
so $\Lambda=\{\lambda\in\Z^n:\ A\lambda=0\}$ is a one-dimensional lattice), then
\[f(\z)=\frac{1-\z^a}{\prod (1-\z^{a_i})},\]
where $a=A\lambda$, and $\lambda$ is the generator of the lattice $\Lambda$.
\end{example}
As in the special case $d=1,n=2$, $\bar{S}$ consists solely of one vertex and one
edge. This formula can also easily be verified directly.

\begin{example}
Theorem \ref{GF}, with $n=d+2$. If the $\R$-span of the $a_i$ is all of $\R^d$ (and
so $\Lambda$ is a two-dimensional lattice), then
\[f(\z)=\sum_j\frac{\z^{p_j}}{(1-\z^{q_j})\prod_i(1-\z^{a_i})}+
\sum_k\frac{\z^{p'_k}}{\prod_i(1-\z^{a_i})},\] where $p_j,q_j,p'_k\in\Z^d$. The
number of terms in the sums is bounded by $C\cdot(nd+\sum \log_2A_{ij})$, for some
constant $C$.
\end{example}
In other words, we can write $f(\z)$ using relatively ``few'' terms.  This is not
immediately obvious, because the number of simplices in $\bar{S}$ may be much larger
than $C\cdot(nd+\sum \log_2A_{ij})$, exponentially larger, in fact. In \cite{S81},
however, H. Scarf shows that $\bar{S}$ has a nice structure, which we will exploit.
In particular, we may represent the edges of $\bar{S}$ by $\{0,h^{ij}\}$, for $i\in
I$ and $0\le j\le N_i$, where $h^{i0},h^{i1},\ldots,h^{iN_i}$ lie on an interval,
that is
\[h^{ij}=c_i+jd_i,\]
for some $c_i,d_i\in\Lambda$. The number of such intervals, $\abs{I}$, is bounded by
$C_1\cdot(nd+\sum \log_2A_{ij})$, where $C_1$ is a constant. The triangles and
3-simplices also lie on intervals (and there are no higher dimensional simplices).
For example, the 3-simplices are
\[\{0,d_i,c_i+(j-1)d_i,c_i+jd_i\},\]
for $i\in I$ and $1\le j\le N_i$.  The exponents in the numerator of
$f_{\bar{S}}(\z)$, which are $A\cdot\max(s)$ for $s\in\bar{S}$, will also lie on
intervals $\alpha_k+j\beta_k$, for $k\in K$, $0\le
 j\le N_k$, and $\alpha_k,\beta_k\in \Z^d$, and we may write
\[\sum_{j=0}^{N_k}\z^{\alpha_k+j\beta_k} \text { as }
\frac{\z^{\alpha_k}-\z^{\alpha_k+(N_k+1)\beta_k}}{1-\z^{\beta_k}}.\] Doing this
gives us a short formula for $f(\z)$.

Here is a specific example:
\begin{example}
Theorem \ref{GF}, with $a_1=(2,0)$, $a_2=(0,3)$, $a_3=(3,8)$, and $a_4=(5,2)$.  Then
\[f(z,w)=\frac{-(z^{20}w^{42}-z^{32}w^6)+(z^{23}w^{50}-z^{35}w^{14})+(z^{22}w^{42}-z^{32}w^{12})-(z^{25}w^{50}-z^{35}w^{20})}
{(1-z^2w^{-6})(1-z^2)(1-w^3)(1-z^3w^8)(1-z^5w^2)}\]
\[+
 \frac{1-z^5w^8-z^{18}w^{48}+z^{20}w^{48}}{(1-z^2)(1-w^3)(1-z^3w^8)(1-z^5w^2)}.\]
 \end{example}
In this example, $\bar{S}$ has one vertex, and it has eight edges on two intervals,
represented by $\{0,h^{ij}\}$, where $h^{10}=(1,-2,1,-1)$ and
\[h^{2j}=(10,14,-5,-1)+(j-1)\cdot(1,-2,1,-1)\text{, for }j=0,\ldots,6.\]
In all, $\bar{S}$ has twelve triangles and five 3-simplices.

Unfortunately, for general $n$ and $d$, the neighborhood complex has no known
structure as nice as in the $n=d+2$ case.  If it did, then perhaps we could write
$f(\z)$ in a short way. For example, L. Lov\'{a}sz conjectured \cite{L89} that the
neighbors of the origin, $b$ such that $\{0,b\}\in S$, are exactly lattice points in
``few'' polytopes of dimension less than $\dim\Lambda$, where ``few'' means the
number is bounded by a polynomial in $nd+\sum \log_2A_{ij}$.  This is the case, as
mentioned, for $n=d+2$, and it is also the case when $n=4,d=1$ (see \cite{S92}), but
for more complicated cases the conjecture is not known to be true or false.

Here is an example of Theorem \ref{GenGF}.
\begin{example}
Theorem \ref{GenGF}, with $a_1=2$, $a_2=3$, and $\Lambda=2 L$, where $
L=\{\ell\in\Z^n:\ A\ell=0\}$ is generated by $(3,-2)$. Then
\begin{align*}
f(\z)&= \frac{1-z^{12}}{(1-z^2)(1-z^3)}\\
&
\begin{array}{l c l c l c l c l c l c l c c}
=&&\phantom{2}1& & &+&\phantom{2}z^2&+&\phantom{2}z^3&+&\phantom{2}z^4&+&\phantom{2}z^5& &\\
&+&2z^6&+&\phantom{2}z^7&+&2z^8&+&2z^9&+&2z^{10}&+&2z^{11}& &\\
&+&2z^{12}&+&2z^{13}&+&2z^{14}&+&2z^{15}&+&2z^{16}&+&2z^{17}&+&\cdots.
\end{array}\\
\end{align*}
\end{example}
In this case, $\Lambda$ is generated by $(6,-4)$, and $\bar{S}$ has one vertex
represented by $\{0\}$ and one edge represented by $\{0,(6,-4)\}$. $T_8$, for
example, contains two points $(4,0)$ and $(1,2)$ (since $8=4\cdot2+0\cdot 3=1\cdot
2+3\cdot 3$). Their difference, $(3,-2)$, is not in $\Lambda$, so $T_8$ has two
equivalence classes, and the coefficient of $z^8$ is 2.  In general, when $d=1$, the
coefficient of $z^a$ is constant for sufficiently large $a$, and it is exactly
$\det(\Lambda)$.  When $d>1$, and if $K\subset\R^d$ is the cone generated by
$a_1,a_2,\ldots,a_n$, the coefficient of $z^a$ is $\det(\Lambda)$ for $a\in K$
sufficiently far from the boundary of $K$.

\section{Proof of Lemma \ref{GF2}}
\begin{figure}
        \begin{center}
        \includegraphics[height= 3in]{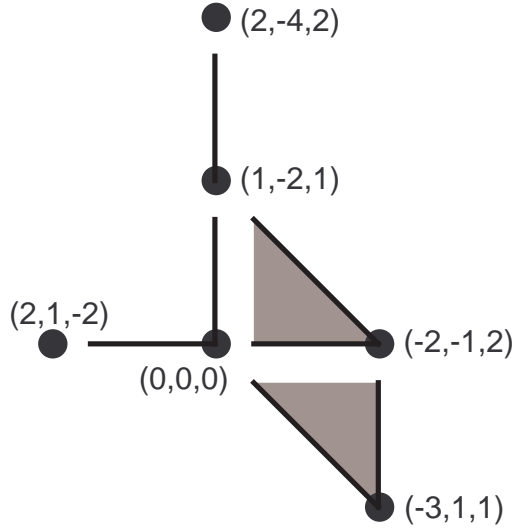}
        \caption{Lemma \ref{Union} applied to Example \ref{ExCXI}}
        \label{FigLemma31}
        \end{center}
\end{figure}

In this section we prove Lemma \ref{GF2}.  Assume that $C$ is a lattice invariant
simplicial complex on $ L=\{\ell\in \mathbb{Z}^n :\  A \ell=0\}$, and let $\bar{C}$
be a set of distinct representatives of $C$ modulo $ L$. We will need the following
basic lemma about $C_\xi$, for $\xi\in\Z^n$, the complex of $s\in C$ such that
$\max(s)\le\xi.$  This lemma says that $C_\xi$ partitions nicely into pieces, and
these pieces are translates of certain subsets of $\bar{C}$.  See Example
\ref{ExLemma31} and Figure \ref{FigLemma31} for an illustration of this lemma
applied to Example \ref{ExCXI}.

\begin{lemma}
\label{Union} Given $\xi\in\Z^n$, and with $ L$ and $C_\xi$ as defined above,
\[C_\xi=\bigcup_{\ell \in  L} \Big( (\bar{C} \cap C_{\xi - \ell}) +
\ell\Big),\] where the union is disjoint.
\end{lemma}
\begin{proof}
Note that the union is disjoint, by the definition of $\bar{C}$.
We will use the fact that
\[C_\xi-\ell=C_{\xi-\ell} \]
for all $\ell \in  L$, since $C$ is invariant under lattice translations. If $s \in
C_\xi$, write $s=s' + \ell$ where $s' \in \bar{C}$ and $\ell \in
 L$.  Then
\[s'=s-\ell \in C_{\xi - \ell}.\]
Therefore $s'\in \bar{C}\cap C_{\xi - \ell}$, and $s\in \Big((\bar{C}\cap C_{\xi -
\ell}) + \ell \Big)$.

Conversely, If $s\in \Big((\bar{C}\cap C_{\xi - \ell}) + \ell \Big)$ for some
$\ell$, then
\[s - \ell \in C_{\xi - \ell},\]
and so $s\in C_\xi$.
\end{proof}

\begin{example}
\label{ExLemma31} Let $C_{\xi}$ be as in Example \ref{ExCXI} and Figure
\ref{FigCXI}. We may take $\bar{C}$ to be the vertex $\{0\}$; the edges
$\{0,(2,1,-2)\}$, $\{0,(1,-2,1)\}$, and $\{0,(3,-1,-1)\}$; and the triangles
$\{0,(2,1,-2),(3,-1,-1)\}$ and $\{0,(1,-2,1),(3,-1,-1)\}$.  Then Lemma \ref{Union}
gives the disjoint union illustrated in Figure \ref{FigLemma31}.
\end{example}

We define another generating function that will be useful in the
proof. Let
\[F_{\bar{C}}(\mathbf{x})=\frac{\sum_{s \in \bar{C}}
(-1)^{\dim(s)}\mathbf{x}^{\max(s)}}{\prod(1-x_i)},
\]
where $\mathbf{x}=(x_1,\ldots,x_n)$.  Then
\[
f_{\bar{C}}(\mathbf{z})=F_{\bar{C}}(\mathbf{z}^{a_1},\mathbf{z}^{a_2},\ldots,
\mathbf{z}^{a_n}).\]

\begin{lemma}
\label{CoeffF} Given $\xi \in \mathbb{Z}^n$ and $F_{\bar{C}}(\x)$ as defined above,
the coefficient of $\mathbf{x}^\xi$ in $F_{\bar{C}}(\mathbf{x})$ is
\[\sum_{s \in \bar{C}\cap C_\xi} (-1)^{\dim(s)}.\]
\end{lemma}
\begin{proof}
For a given $s\in \bar{C}$, the term
\[\frac{(-1)^{\dim(s)}\mathbf{x}^{\max(s)}}{\prod(1-x_i)}
\]
will contribute $(-1)^{\dim(s)}\mathbf{x}^\xi$ if $\max(s)\le
\xi$, and otherwise it will contribute nothing.  The proof
follows, by the definition of $C_\xi$.
\end{proof}

Now we have the tools to prove Lemma \ref{GF2}.

\

{\bf Proof of Lemma \ref{GF2}:\ } Given a $ L$-invariant simplicial complex, $C$,
fix $b\in\Z^d$.  Take a particular $\xi_0\in\Z^n$ such that $A\xi_0 = b$. Then all
$\xi\in \mathbb{Z}^n$ such that $A\xi =b$ are given by $\xi_0 - \ell$, for $\ell \in
L$. Let $d_b=\text{EC}(C_{\xi_0}).$ We want to show that the coefficient of $\z^b$
in $f_{\bar{C}}(\z)$ is $d_b$.  Since
$f_{\bar{C}}(\mathbf{z})=F_{\bar{C}}(\mathbf{z}^{a_1},\mathbf{z}^{a_2},\ldots,
\mathbf{z}^{a_n})$,
\begin{align*}
    \text{the coefficient of }\mathbf{z}^b \text{ in } f_{\bar{C}}(\mathbf{z}) & =
\sum_{\xi\in\Z^n:\ A\xi=b}\text{the coefficient of }\x^\xi \text{ in }
F_{\bar{C}}(\x)\\
    &=\sum_{\ell \in \Lambda} \text{the coefficient of }\mathbf{x}^{\xi_0-\ell} \text{ in }
    F_{\bar{C}}(\mathbf{x})\\
\end{align*}
\begin{align*}
    \phantom{\text{the coefficient of }\mathbf{z}^b \text{ in } f_{\bar{C}}(\mathbf{z})}& = \sum_{\ell \in \Lambda}\sum_{\genfrac{}{}{0pt}{}{s'\in}{\bar{C}\cap
    C_{\xi_0-\ell}}}(-1)^{\dim(s')}\text{    (by Lemma \ref{CoeffF})}\\
    &=\sum_{\genfrac{}{}{0pt}{}{s \in C_{\xi_0},}{\ell:\ s-\ell\in\bar{C}}}(-1)^{\dim(s-\ell)}\text{   (by Lemma \ref{Union})}\\
    &=\sum_{s \in C_{\xi_0}}(-1)^{\dim(s)}\\
    &=d_b.\\
\end{align*}

\

We have proven that, for all $b$, the coefficient of $z^b$ is the same in $\sum
d_b\z^b$ and in $f_{\bar{C}}(\mathbf{z})$, and the proof follows.
\hspace{\stretch{1}}$\Box$

\section{The Neighborhood Complex}
Assume that $\Lambda$ is a generic lattice such that $A\lambda=0$ for all
$\lambda\in\Lambda$ (we will deal with the nongeneric case in Section 5), and let
$S$ be the neighborhood complex, as defined in Section 1. In this section, we will
prove Theorems \ref{GF} and \ref{GenGF}. First we will examine $S$ and the
subcomplexes $S_\xi$ (the complex of $s\in S$ such that $\max(s)\le\xi$). Our goal
is to prove the following lemma.

\begin{lemma}
\label{ECS} Given $S$ as above, for $\xi\in\Z^n$, if $S_\xi\ne \emptyset$, then
EC$(S_\xi)=1$.
\end{lemma}

We will prove this lemma by giving a geometric realization of the $S_\xi$ and then
using properties of this realization to compute the Euler characteristic. We will
use a construction from \cite{BSS98}, where the authors prove that a particular
complex (the neighborhood complex with ideal vertices included) is homeomorphic to
$\mathbb{R}^{m-1},$ where $m=\dim(\Lambda)$. In fact, the $S_\xi$ also have a nice
topological property: they are contractible (this is shown in \cite{BPS}).
Contractibility implies that the Euler characteristic is 1 (this can be seen by
applying standard facts from the homology of CW-complexes, see, for example, Theorem
IX.4.4 of \cite{M91}), but here we will find EC($S_\xi$) directly and geometrically.
Bayer and Sturmfels \cite{BS98} also use a very similar construction to analyze
their \emph{hull complex}.

For purposes of exposition, we will present lemmas in a different order from how
they are proved. The structure of the proof of Lemma \ref{ECS} is: Lemma
\ref{LemmaBoundedFaces} and Lemma \ref{Polyhedron} imply Lemma \ref{EC}, and then
Lemma \ref{CX} and Lemma \ref{EC} imply Lemma \ref{ECS}.

Let $X=\{x^1,x^2,\ldots,x^m\}$, with $x^i\in\mathbb{R}^n$, be given.  We define the
complex $C(X)$ on the vertices $X$ to be the $s\subset X$ such that there is no
$x\in X$ with $x<\max(s)$. $C(X)$ is a simplicial complex. We first prove the
following lemma.

\begin{lemma}
\label{CX} For $\xi\in\Z^n$, if $X=\{x\in\Lambda:\  x\le \xi\}$, then $S_\xi=C(X)$.
\end{lemma}
\begin{proof}
Suppose $s=\{\lambda^1,\lambda^2,\ldots,\lambda^k\}\in S_\xi$.
Then $\lambda_1,\ldots,\lambda_k\le \xi$ and for no
$\lambda\in\Lambda$ is $\lambda<\max(s)$.  Therefore for no $x\in
X$ is $x<\max(s)$ (since $X\subset\Lambda$), and so $s\in C(X)$.

Conversely, suppose $s=\{\lambda^1,\lambda^2,\ldots,\lambda^k\}\in
C(X)$.  Then $\lambda_1,\ldots,\lambda_k\le \xi$ and for no $x\in
X$ is $x<\max(s)$.  Suppose (seeking a contradiction) that
$\lambda < \max(s)$ for some $\lambda\in\Lambda$.  Then for each
$i$ there is a $j$ such that
\[\lambda_i<\lambda^j_i\le \xi_i.\]
But then $\lambda < \xi$ and so $\lambda\in X$, contradicting that
for no $x\in X$ is $x<\max(s)$.  Therefore, for no
$\lambda\in\Lambda$ is $\lambda<\max(s)$, and so $s\in S_\xi$.
\end{proof}

We say that $X$ is \emph{generic} if, whenever there is some $x^1,x^2\in X$, with
$x^1\ne x^2$ but $x^1_i=x^2_i$ for some $i$, then there is an $x\in X$ with
$x<\max(x^1,x^2)$.  This definition is slightly more complicated than for a lattice,
because $X$ need not be lattice invariant.  Then Lemma \ref{ECS} will follow from
Lemma \ref{CX}, and the following lemma.

\begin{lemma}
\label{EC} If $X=\{x^1,x^2,\cdots,x^m\}$ is generic and $C(X)$ is
defined as above, then EC$\big(C(X)\big)=1$.
\end{lemma}

To prove this lemma, we follow the method of \cite{BSS98} and construct a polyhedron
$P_t$ from the points $x^1,x^2,\ldots,x^m$, as follows. Given $t\ge 0$, define
$E_t:\ \mathbb{R}^n\rightarrow\mathbb{R}^n$ by
\[E_t(x)=\mathbf{e}^{tx}=(e^{tx_1},e^{tx_2},\ldots,e^{tx_n}),\]
where $x=(x_1,x_2,\ldots,x_n)$. Now we define
\[P_t=\mathbb{R}^n_{\ge 0} +\text{conv}\{E_t(x^1),E_t(x^2),\ldots,E_t(x^m)\},\]
where $X=\{x^1,x^2,\ldots,x^m\}$.

\begin{figure}
        \begin{center}
        \includegraphics[height= 3in]{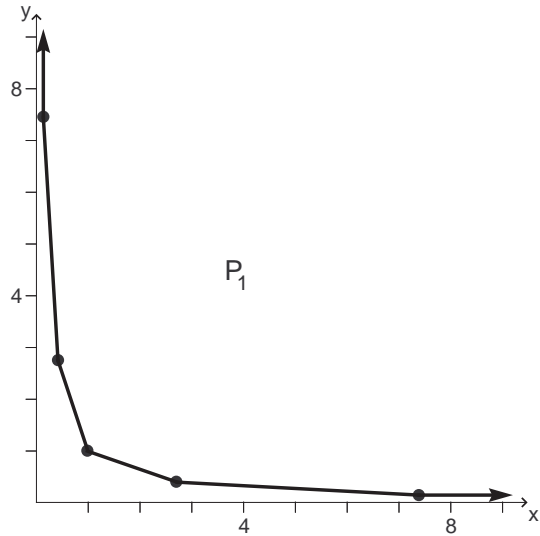}
        \caption{Example \ref{Ex3}, $P_1$, with $X=\{-2,-1,0,1,2\}$}
        \label{Fig3}
        \end{center}
\end{figure}

\begin{example}\label{Ex3} Let $X=\{-2,-1,0,1,2\}$.  Then
Figure \ref{Fig3} illustrates $P_1$.
\end{example}

The polyhedron $P_t$ has the following useful property.

\begin{lemma}
\label{Polyhedron} There exists a sufficiently large $t$ such that, if $s\subset X$
with $s=\{s^1,s^2,\ldots,s^k\}$, then $s\in C(X)$ if and only if
conv$\{E_t(s^1),E_t(s^2),\ldots,E_t(s^k)\}$ is a face of $P_t$.
\end{lemma}

\begin{proof}
The proof is very similar to the proof of Theorem 2 of \cite{BSS98}. We won't go
through the details.
\end{proof}

In Example \ref{Ex3} (see Figure \ref{Fig3}), this lemma tells us that $C(X)$ has
vertices $-2,-1,0,1,2$, and edges $\{-2,-1\},\{-1,0\},\{0,1\},\{1,2\}$, as we would
expect.  In general, Lemma \ref{Polyhedron} gives a geometric realization of $C(X)$
in $\R^{n}$.  In fact, as shown in Theorem 2 of \cite{BSS98}, if we take $X$ to be
the (infinite) set $\Lambda$, $P_t$ gives a geometric realization of $S$, the entire
neighborhood complex.

Now pick a sufficiently large $t$ such that Lemma \ref{Polyhedron} holds.  Then the
simplices in $C(X)$ are exactly the bounded faces of $P_t$.  Then Lemma \ref{EC}
(and hence Lemma \ref{ECS}) follows from the following lemma.

\begin{lemma}\label{LemmaBoundedFaces}  Let $P$ be an unbounded polyhedron in $\mathbb{R}^n$. Let $\mathcal{F}$ be the
collection of bounded faces of $P$. Then
\[\text{EC}(\mathcal{F})=\sum_{F\in\mathcal{F}}(-1)^{\dim(F)}=1.\]
\end{lemma}

\begin{proof}
Choose a half-space $H_+$ such that $H_+$ contains all of the bounded faces of $P$
in its interior and such that $P'=P\cap H_+$ is bounded. Let $\mathcal{F}'$ be the
collection of faces of $P'$. We know
\[\sum_{F'\in\mathcal{F}'}(-1)^{\dim(F')}=1+(-1)^{n-1}.\]
This is the Euler-Poincar\'{e} formula, and it can be seen
combinatorially (see, for example, Corollary VI.3.2 of
\cite{B02}), or it can be seen from the fact that the complex
$\mathcal{F}'$ is homeomorphic to an $n-1$ sphere (and then
applying standard facts from the homology of CW-complexes, see,
for example, Theorem IX.4.4 of \cite{M91}). Let $H$ be the
hyperplane which is the boundary of $H_+$. The faces of $P'$ fall
into $4$ categories:
\begin{enumerate}
    \item $\mathcal{F}$, the bounded faces of $P$,\\
    \item The face $P\cap H$,\\
    \item $F\cap H_+$, where $F$ is an unbounded face of $P$, and\\
    \item $F\cap H$, where $F$ is an unbounded face of $P$.\\
\end{enumerate}
There is a bijective correspondence between the last two categories, mapping a face
$F$ from category $3$ of dimension $k$ to $F\cap H$, a face from category $4$ of
dimension $k-1$. Therefore, in $\sum_{F'\in\mathcal{F}'}(-1)^{\dim(F')}$, these two
categories will exactly cancel each other, and so we have
\[1+(-1)^{n-1}=\sum_{F'\in\mathcal{F}'}(-1)^{\dim(F')}=\Big[\sum_{F\in\mathcal{F}}(-1)^{\dim(F)}\Big]+(-1)^{n-1}+0.\]
The lemma follows.
\end{proof}

Now we are ready to prove Theorems \ref{GF} and \ref{GenGF} (in the generic case).

\

{\bf Proof of Theorem \ref{GF}:} Let $\Lambda= L=\{\lambda\in\Z^d:\ A\lambda=0\}$,
and let $S$ be the neighborhood complex on $\Lambda$.  Take a particular
$\xi_0\in\Z^n$ such that $A\xi_0=b$, and let $d_b=\text{EC}(S_{\xi_0})$.  We want to
show that
\[f_{\bar{S}}(\z)=\sum_{b\in G}\z^b,\]
and by Lemma \ref{GF2} we know that
\[f_{\bar{S}}(\z)=\sum_{b\in \Z^d} d_b\z^b.\]
By Lemma \ref{ECS}, we know that $d_b=1$ if and only if $S_{\xi_0}$ is nonempty (and
$d_b=0$ otherwise), so it suffices to show that $S_{\xi_0}$ is nonempty if and only
if $b\in G$.

Indeed, if $\{\lambda\}\in S_{\xi_0}$, for some $\lambda\in\Lambda$, then
$\lambda\le\xi_0$ and so $\xi_0-\lambda\ge 0$.  Then, since
$A\cdot(\xi_0-\lambda)=b-0=b$ with $\xi_0-\lambda\ge 0$, we have that $b\in G$.
Conversely, if $b\in G$, then there is some $\xi\ge 0$ such that $A\xi=b$.  Then
$\xi_0-\xi\le \xi_0$, and $A\cdot(\xi_0-\xi)=b-b=0$, so $\xi_0-\xi\in\Lambda$ and
$\{\xi_0-\xi\}\in S_{\xi_0}$.  The proof follows. \hspace{\stretch{1}}$\Box$

\

{\bf Proof of Theorem \ref{GenGF}:} Let $\Lambda$ be a lattice in $\Z^n$ such that
$A\lambda=0$, for all $\lambda\in\Lambda$, and let $S$ be the neighborhood complex
defined on $\Lambda$. Recall that, for $b\in \Z^n$, we define $T_b=\{\xi\in\Z^n:\
\xi\ge 0 \text{ and } A\xi=b\}$, we define an equivalence relation on $T_b$ by
$\xi\sim\eta$ if and only if $\xi-\eta\in\Lambda$, and we define $c_b$ to be the
number of equivalence classes in $T_b$.  To use Lemma \ref{GF2}, we must have a
lattice invariant neighborhood complex on all of $ L=\{\ell\in \mathbb{Z}^n :\  A
\ell=0\}$.  Let $\bar{L}$ be a set of distinct representatives of $L$ modulo
$\Lambda$, and define $C$ to be the disjoint union
\[C=\bigcup_{\ell\in\bar{L}}S+\ell.\]
$C$ is an $ L$-invariant complex, and we can choose $\bar{C}$ and $\bar{S}$
(representatives of $C$ modulo $ L$ and $S$ modulo $\Lambda$, respectively) such
that $\bar{C}=\bar{S}$. By Lemma \ref{GF2}, we know
\[f_{\bar{S}}(\z)=f_{\bar{C}}(\z)=\sum_{b\in\Z^d}d_b\z^b,\]
where $d_b=\text{EC}(C_\xi),\text{ for some (any) }\xi\text{ such that } A\xi=b$.
Therefore we need to show that $c_b=d_b$, for all $b\in\Z^n$.

Fix a $\xi_0$ such that $A\xi_0=b$.  We claim that
\[C_{\xi_0}=\bigcup_{\ell\in\bar{L}}(S_{\xi_0+\ell}-\ell),\]
where the union is disjoint.  Indeed, if $s=\{\ell^1,\ell^2,\ldots,\ell^k\}\in
C_{\xi_0}$, then, for each $i$, $\ell^i\le\xi_0$. Take $\ell\in \bar{L}$ such that
$s+\ell\in S$. Then $\ell^i+\ell\le\xi_0+\ell$, and so $s+\ell\in S_{\xi_0+\ell}$
and $s\in S_{\xi_0+\ell}-\ell.$ Conversely, if $s=\{\ell^1,\ell^2,\ldots,\ell^k\}\in
S_{\xi_0+\ell}-\ell,$ then $\ell^i+\ell\le\xi_0+\ell$, for all $i$. Therefore,
$\ell^i\le\xi_0$, and so $s\in C_{\xi_0}$.  In addition, the union is disjoint,
because $S_{\xi+\ell}+\ell\subset S+\ell$, which are themselves disjoint.

Since we have written $C_{\xi_0}$ as a disjoint union, we have
\[\text{EC}(C_{\xi_0})=
\sum_{\ell\in\bar{L}}\text{EC}(S_{\xi_0+\ell}).\] Since EC$(S_{\xi_0+\ell})=1$ if
$S_{\xi_0+\ell}\ne \emptyset$, by Lemma \ref{ECS}, and EC$(S_{\xi_0+\ell})=0$ if
$S_{\xi_0+\ell}= \emptyset$, we have
\[\text{EC}(C_{\xi_0})=\#\{\ell\in\bar{L}:\ S_{\xi_0+\ell}\ne
\emptyset\}.\] Therefore, to prove Theorem \ref{GenGF}, we must show that the number
of nonempty $S_{\xi_0+\ell}$ is the number, $c_b$, of equivalence classes of $T_b$.

For $\xi,\eta\in T_b$, $\xi\sim\eta$ if and only if $\xi-\eta\in\Lambda$, which
happens if and only if $\xi$ and $\eta$ are in the same coset
$\xi_0+(\ell+\Lambda)$, for some $\ell\in\bar{L}$. Then the equivalence classes of
$T_b$ are exactly the $\big(\xi_0+(\ell+\Lambda)\big)\cap\Z^n_{\ge 0}$ which are
nonempty. But $\lambda\in\Lambda$ is such that $\xi_0+(\ell+\lambda)\ge 0$ if and
only if $-\lambda\le \xi_0+\ell$, which happens if and only if $\{-\lambda\}\in
S_{\xi_0+\ell}$.  Therefore $S_{\xi_0+\ell}$ is nonempty if and only if
$\big(\xi_0+(\ell+\Lambda)\big)\cap\Z^n_{\ge 0}$ is a nonempty equivalence class of
$T_b$.  The proof of Theorem \ref{GenGF} follows. \hspace{\stretch{1}}$\Box$

\section{The Non-generic Case}
The strategy we follow is to perturb the elements of $\Lambda$ so
that no two have any coordinate that is the same.  Then we will be
in the generic case and can apply the lemmas of the last section.

We call $\varphi:\Lambda\rightarrow\mathbb{R}^n$ a \emph{proper
perturbation} if the following 3 conditions hold:
\begin{enumerate}
   \item If $x\ne y$, then $[\varphi(x)]_i \ne [\varphi(y)]_i$,
   \item If $[\varphi(x)]_i<[\varphi(y)]_i$, then $x_i\le y_i$, and
   \item If $[\varphi(x)]_i<[\varphi(y)]_i$, then
   $[\varphi(x+\lambda)]_i<[\varphi(y+\lambda)]_i$ for all
   $\lambda\in\Lambda$.
\end{enumerate}

The first condition insures that we will be in the generic case, the second insures
that the perturbation only ``breaks ties'' and doesn't change the natural ordering,
and the third condition will be needed to prove that the neighborhood complex is
lattice invariant.

To prove that proper perturbations exist, we will construct an example of one.

\begin{example} This example corresponds to the lexicographical
tie-breaking rule used in \cite{S81}. Given an integer $i$, let
$f_i:\mathbb{Z}\rightarrow\mathbb{R}$ be a function such that
\begin{enumerate}
    \item $f_i$ is strictly increasing,
    \item $f_i(0)=0$ (an hence $f_i(x)<0$ if $x<0$), and
    \item if $\abs{x}>0$ (hence $\abs{x}\ge 1$), then $\frac{1}{2^{2i}}\le \abs{f_i(x)} <
    \frac{1}{2^{2i-1}}$.
\end{enumerate}
For example, $f_i$ could be an appropriate rescaling of $\arctan(x)$.  Now define
$\varphi:\Lambda\rightarrow\mathbb{R}^n$ by
\[\varphi(x)=x+(x_1f_1(x_1)+x_2f_2(x_2)+\cdots+x_nf_n(x_n))\cdot\mathbf{1},\]
where $\mathbf{1}$ is the $n$-vector of ones. One can check that $\varphi$ is a
proper perturbation.
\end{example}

Given a proper perturbation $\varphi$, we can now define the neighborhood complex,
$S$, on the vertices $\Lambda$, by saying
$s=\{\lambda^1,\lambda^2,\ldots,\lambda^k\}$ is in $S$ if and only if for no
$\lambda \in \Lambda$ is $\varphi(\lambda)<\max(\varphi(s))$, where
$\varphi(s)=\{\varphi(\lambda^1),\varphi(\lambda^2),\ldots,\varphi(\lambda^k)\}$.
$S$ may be different for different $\varphi$, but many properties (including
Theorems \ref{GF} and \ref{GenGF}) hold regardless of the choice of $\varphi$. The
following lemma shows that $S$ is invariant under lattice translations, and so
$f_{\bar{S}}(\z)$, as defined in Section 1, makes sense.

\begin{lemma}
If $\varphi$ is a proper perturbation, then the neighborhood complex $S$, as defined
above, is lattice invariant.
\end{lemma}
\begin{proof}
Given $\lambda\in\Lambda$, we have the following chain of
implications:
\begin{align*}
    s&=\{\lambda^1,\lambda^2,\ldots,\lambda^k\}\in S\\
    &\Rightarrow\text{for no }\lambda'\in\Lambda\text{ is }
    \varphi(\lambda')<\max(\varphi(s))\\
&\Rightarrow \text{given }\lambda'\in\Lambda, \ \exists i\text{ such that } \forall j\ [\varphi(\lambda')]_i\ge [\varphi(\lambda^j)]_i\\
&\Rightarrow \text{given }\lambda'\in\Lambda, \ \exists i\text{ such that } \forall j\ [\varphi(\lambda'+\lambda)]_i\ge [\varphi(\lambda^j+\lambda)]_i\\
&\phantom{XXX}\text{(by Property 3 of proper perturbations)}\\
&\Rightarrow\text{for no }\lambda'\in\Lambda\text{ is }
    \varphi(\lambda'+\lambda)<\max(\varphi(s+\lambda))\\
&\Rightarrow s+\lambda\in S.
\end{align*}
\end{proof}

Given $\xi\in\mathbb{Z}^n$, we define $S_\xi$ as in Section 1, that is, $S_\xi$ is
the complex of all $s\in S$ such that $\max(s)\le\xi$. For generic $X\subset \R^n$,
define $C(X)$ as in Section 4, that is, $C(X)$ is the simplicial complex of
$s\subset X$ such that there is no $x\in X$ with $x<\max(s)$. We mimic Lemma
\ref{CX}.

\begin{lemma}
If $\varphi$ is a proper perturbation, if $\xi\in\Z^n$ is given, and if
$Y=\{y\in\Lambda:\  y\le \xi\}$, then $\varphi(S_\xi)=C(\varphi(Y))$ (and hence
$S_\xi$ is isomorphic to $C(\varphi(Y))$).
\end{lemma}
\begin{proof}
Suppose $s=\{\lambda^1,\lambda^2,\ldots,\lambda^k\}\in S_\xi$. Then
$\lambda^1,\ldots,\lambda^k\le \xi$ and for no $\lambda\in\Lambda$ is
$\varphi(\lambda)<\max(\varphi(s))$. Therefore for no $y\in Y$ is
$\varphi(y)<\max(\varphi(s))$ (since $Y\subset\Lambda$), and so $\varphi(s)\in
C(\varphi(Y))$.

Conversely, suppose $\varphi(s)\in C(\varphi(Y))$, with
$s=\{\lambda^1,\lambda^2,\ldots,\lambda^k\}$.  Then $\lambda^1,\ldots,\lambda^k\le
\xi$ and for no $y\in Y$ is $\varphi(y)<\max(\varphi(s))$.  Suppose (seeking a
contradiction) that $\varphi(\lambda) < \max(\varphi(s))$ for some
$\lambda\in\Lambda$.  Then for each $i$ there is a $j$ such that
\[[\varphi(\lambda)]_i<[\varphi(\lambda^j)]_i.\]
Therefore $\lambda_i\le\lambda^j_i$, by Property 2 of proper
perturbations, and so
\[\lambda_i\le\lambda^j_i\le \xi_i.\]
But then $\lambda \le \xi$ and so $\lambda\in Y$, contradicting that for no $y\in Y$
is $\varphi(y)<\max(\varphi(s))$. Therefore, for no $\lambda\in\Lambda$ is
$\varphi(\lambda)<\max(\varphi(s))$, and so $s\in S_\xi$.
\end{proof}

In particular, this lemma, together with Lemma \ref{EC}, implies that
\[\text{EC}(S_\xi)=\text{EC}\big(C(\varphi(Y))\big)=1\]
whenever $S_\xi$ is nonempty.  The proofs of Theorems \ref{GF} and \ref{GenGF} in
the non-generic case are now identical to their proofs in the generic case (see
Section 4).

\ \\

\noindent \small Cowles Foundation for Research in Economics, Yale University, New Haven, Connecticut 06511\\
herbert.scarf@yale.edu\\

\noindent Department of Mathematics, University of Michigan, Ann Arbor, Michigan
48109\\kmwoods@umich.edu\\

\end{document}